\documentclass[12pt, leqno]{elsart}

\usepackage[mathcal]{euscript}
\usepackage{amssymb, amsfonts, xypic}
\usepackage{graphicx}

\xyoption{all}
\CompileMatrices

\begin{document}
\begin{frontmatter}

\title{Homotopical localizations at a space\thanksref{part_of_thesis}}
\thanks[part_of_thesis]{
This paper forms a modified version of a part of the
author's Ph.D. dissertation written at the University of Chicago
under the direction
of J. P. May; current address
of the author: Warsaw University}
\author{Adam J. Prze\'zdziecki}
\address{
Department of Mathematics, University of Chicago, Chicago, USA \\
Institute of Mathematics, Warsaw University, Warsaw, Poland}
\ead{adamp@mimuw.edu.pl}

\begin{abstract}
Our main motivation for the work presented in this paper is to
construct a localization functor, in a certain sense dual to the
$f$-localization of Bousfield and Farjoun, and to study some of its
properties. We succeed in a case which is related to the Sullivan
profinite completion. As a corollary we prove the existence of
certain cohomological localizations.

\vspace{4pt}
{\noindent\leavevmode\hbox {\it AMS classification:\ }}
55P60; 55P65
\vspace{-9pt}
\end{abstract}

\begin{keyword}
homotopical localization \sep cohomological localization \sep
$f$-localization
\end{keyword}

\pagestyle{plain}
\lefthyphenmin 4
\righthyphenmin 4

\end{frontmatter}

\theoremstyle{definition}
\newtheorem{definition}[thm]{Definition}
\newtheorem{example}[thm]{Example}

\theoremstyle{remark}

\newcommand{\map}{\mbox{map}}
\newcommand{\colim}{\mbox{colim}}
\newcommand{\holim}{\mbox{holim}}
\newcommand{\hocolim}{\mbox{hocolim}}
\newcommand{\comment}[1]{}

\makeatletter
\let\c@equation\c@thm
\makeatother

\begin{center}
\end{center}
\maketitle

\section{Introduction}

We can view $f$-localization as the initial coaugmented idempotent
functor on the homotopy category which takes a map
$f$ to an equivalence. In
\cite{bousfield_factorization} Bousfield used the small object
argument to prove that $f$-localizations exist for all maps $f$.
The role of
these functors was especially exposed in $1990$'s when they were put
in a convenient framework in terms of mapping complexes.
A survey of related methods can be found in \cite{farjoun} and
\cite{casacuberta_recent}. It seems
natural to ask if a dual notion of a localization at a space $Z$, that
is the terminal idempotent functor with a given space $Z$ in its image
(Definition \ref{def_localization_at_z}), might not also be
interesting. The main reason these localizations
have not been considered very much
is that they are not known to exist in general, even in the stable
case (see Chapter $7$ in \cite{margolis}).

As every homological localization can be realized as an
$f$-localization, every cohomological localization, provided it exists,
is a localization at a suitable space. Research towards establishing
the existence of cohomological localizations was briefly summarized in
$2.6$ of \cite{bousfield_homotopical_spaces}.

Here we prove the existence of
localizations at compactly topologized spaces
(Definition \ref{definition_compact} and
Theorem \ref{theorem_localization_at_Z}).
Examples of such spaces
include the ones which are profinite completions of another space,
mapping complexes with a profinitely completed target, and others.
This result allows us to construct an idempotent approximation to the
Sullivan profinite completion
(Theorem \ref{approximation_to_completion}).

We would like to be able to prove the existence of localization at
an arbitrary space without relying on the compactness condition,
and there is some evidence that such localizations should exist at abelian
Eilenberg-Mac Lane spaces. These would form
``truncated localizations at an ordinary cohomology theory'', an
analogue of ``truncated localizations at a homology theory''
whose existence was shown by Ohkawa in \cite{ohkawa}.
It would also be interesting to find how such localizations act on
spaces and how they are related to those $f$-localizations,
that do not correspond to a localization at any space.

Casacuberta, Scevenels and Smith
investigated in \cite{casacuberta_cardinal}
dependence on certain large cardinal axioms of
a more general question, from a positive answer to which
the existence of localizations at any space would follow.
Despite extensive efforts we were unable to avoid similar set theoretic
problems in our attempts to
prove the existence of localizations at a general space, nor
were we able to disprove it under some large cardinal axioms.

The main Theorem
\ref{theorem_localization_at_Z} is proved in section
\ref{section_at_a_space}. In section
\ref{section_applications_and_examples} we describe an idempotent
approximation to the Sullivan profinite completion and prove
the existence of certain cohomological localizations.

The paper is written simplicially.
We use terms ``space'' and ``simplicial
set'' as synonyms choosing the second one wherever confusion with
compact topological space might occur or to emphasize it when we work
on the point set level rather than in the homotopy category.
To make the presentation more accessible, we frequently work in
the pointed homotopy category $Ho_*$. Adjective ``compact''
always means ``compact Hausdorff''.

The author is very grateful to J. P. May and A. K. Bousfield for
their support, encouragement and helpful comments.

\section{Localizations} \label{sec_localizations}

In this section we collect basic definitions and facts related to
homotopical localizations.

A functor $L$ is called {\em coaugmented} if it comes with a natural
transformation $\eta_X:X\to LX$ from the identity to $L$. A coaugmented
functor is {\em idempotent} if
in the diagram
$$
\xymatrix{
  X  \ar[d] \ar[r] & LX \ar[d]^{\eta_{LX}} \\
  LX \ar[r]^{L\eta_X} & LLX
  }
$$
the maps $\eta_{LX}$ and $L\eta_X$ are equivalences and
$\eta_{LX}=L\eta_X$.
\begin{defn}
  A coaugmented idempotent functor is called a {\em localization}.
\end{defn}

Although this definition makes sense in any category we will consider
only localizations in the homotopy category
$Ho_*$ of pointed simplicial sets (spaces).
A space $Z$ is said to be {\em $L$-local} if the map
$\eta_Z:Z\to LZ$ is an equivalence. It is straightforward to check that
the class of $L$-local spaces uniquely determines and is determined by the
functor $L$.
A map $g:X\to Y$ is an {\em $L$-equivalence} if $Lg$ is an
equivalence. There is a natural ordering of localizations as described
below.
\begin{defn}
  Given two localization functors $L_1$ and $L_2$ we say that
  $L_1\leq L_2$ if one of the equivalent conditions hold:
  \begin{itemize}
    \item[(i)] there is a natural transformation $L_1\to L_2$ giving
      $L_2L_1\simeq L_2$
    \item[(ii)] any $L_1$-equivalence is also an $L_2$-equivalence
    \item[(iii)] any $L_2$-local space is also $L_1$-local
  \end{itemize}
\end{defn}
This definition is an obvious extension of the ordering in the Bousfield
lattice of $f$-localizations ($4.3$ in
\cite{bousfield_homotopical_spaces}).

Given a map $\nolinebreak{f:A\to B}$ we say that
a fibrant space $Z$ is
{\em $f$-local} if the induced map of function complexes
\begin{equation}\label{equation_f_local}
f^*:\map_*(B,Z)\to\map_*(A,Z)
\end{equation}
is an equivalence. If $Z$ is connected the condition above is
equivalent to the one that the induced map of unbased function
complexes
$$
f^*:\map(B,Z)\to\map(A,Z)
$$
is an equivalence.

A map $g:X\to Y$ is an
$\nolinebreak{f\mbox{\em{}-equivalence}}$ if any $f$-local space is
also $g$-local. This means that for any fibrant space $Z$ if
$$
f^*:\map_*(B,Z)\stackrel{\simeq}{\to}\map_*(A,Z)
$$
then
$$
g^*:\map_*(Y,Z)\stackrel{\simeq}{\to}\map_*(X,Z).
$$
\begin{defn}
  An {\em $f$-localization} is a localization functor $L_f$ such that the
  following conditions hold:
  \begin{itemize}
  \item[(i)] The classes of $f$-equivalences and $L_f$-equivalences
    coincide.
  \item[(ii)] The classes of $f$-local and $L_f$-local spaces
    coincide.
  \item[(iii)] The map $X\to L_fX$ is an $f$-equivalence and $L_fX$ is
    $f$-local.
  \item[(iv)] $L_f$ is the initial localization functor such that the
    map $f$ is an $L_f$-equivalence.
  \end{itemize}
  For a map $f$, there are obvious implications
  $(i) \Leftrightarrow (ii) \Leftrightarrow (iii) \Leftrightarrow (iv)$.
\end{defn}
The existence of $f$-localizations for arbitrary maps $f$ was proved
by Bousfield \cite{bousfield_factorization} and Farjoun
\cite{farjoun_map}.

Let $Z$ be a fibrant space. We say that a map $g:X\to Y$ is a
$\nolinebreak{Z\mbox{\em{}-equivalence}}$
if the induced map of function complexes
$$
g^*:\map_*(Y,Z)\to\map_*(X,Z)
$$
is an equivalence. A fibrant space $K$ is {\em $Z$-local} if it is
$g$-local for all $Z$-equivalences $g$. This means that for any $g$ if
$$
g^*:\map_*(Y,Z)\stackrel{\simeq}{\to}\map_*(X,Z)
$$
then
$$
g^*:\map_*(Y,K)\stackrel{\simeq}{\to}\map_*(X,K)
$$
\begin{defn} \label{def_localization_at_z}
  A {\em localization at} $Z$ is a localization functor $L_Z$ such that the
  following conditions hold:
  \begin{itemize}
  \item[(i)]  the classes of $Z$-equivalences and $L_Z$-equivalences
    coincide.
  \item[(ii)] the classes of $Z$-local and $L_Z$-local spaces coincide.
  \item[(iii)] The map $X\to L_ZX$ is a $Z$-equivalence and $L_ZX$ is
    $Z$-local.
  \item[(iv)] $L_Z$ is the terminal localization functor such that
      the space $Z$ is $L_Z$-local.
  \end{itemize}
  For a space $Z$, there are obvious implications
  $(i) \Leftrightarrow (ii) \Leftrightarrow (iii) \Rightarrow (iv)$.
\end{defn}
The implication $(iv) \Rightarrow (iii)$ is obvious when $L_Z$
in the sense of {\em (i) - (iii)} exists.
The only problem might arise if $L_Z$ exists in the sense of
{\em (iv)} but not {\em (i) - (iii)}, that is, a terminal
localization $T$ such that $Z$ is $T$-local exists but not
all $T$-local spaces are $Z$-local (condition {\em (ii)}).
Suppose $K$ is such a $T$-local but not $Z$-local space.
Then there is a $Z$-equivalence $f:A\to B$ which is not a $K$-equivalence.
Thus $K$ is $T$-local but not $f$-local hence $L_f$ is not
less than $T$ which contradicts {\em (iv)}.

The existence of localization at a given space $Z$ is not known in
general.

\comment{The work of Casacuberta, Scevenels and Smith
\cite{casacuberta_cardinal} implies that the existence of certain
homotopical localizations may depend on large cardinal axioms.

that under the large cardinal axiom called Vopenka principle this
ordering is naturally isomorphic to the Bousfield lattice of
$f$-localizations ($4.3$ in \cite{bousfield_homotopical_spaces}).
}

It is clear that the classes of $Z$-equivalences and
$f$-equivalences are closed under arbitrary homotopy colimits.
Also the classes of $Z$-local and $f$-local spaces are closed under
arbitrary homotopy limits.
\begin{lem}\label{lemma_set_generated_equivalences}
  Suppose that for a certain space $Z$ there is a set of
  $Z$-equivalences $\{f_\alpha\}$ such that every $Z$-equivalence can be
  presented as a homotopy colimit of elements of the set $\{f_\alpha\}$.
  Then the localization at $Z$ is
  simply an $f$-localization for $f=\bigvee f_\alpha$.
\end{lem}

\section{A characterization of $Z$-equivalences} \label{sec_characterisation}

In this section we recall Lemma \ref{lemma_llp}.
Although it is not new we prove it here since we didn't find
an appropriate reference.

We say that a map $f:A\to B$ has a left lifting property (LLP)
with respect to a map $g:C\to D$ if any diagram
$$
\xymatrix{
  A \ar[r] \ar[d]_f & C \ar[d]^g \\
  B \ar@{-->}[ur] \ar[r] & D
}
$$
admits the dashed map. For the sake of clarity we will use the term
homotopy LLP when the lift we have in mind is in the homotopy
category.

\begin{lem} \label{lemma_symmetry}
  Let $f:A\to B$ and $g:C\to D$ be maps in $Ho_*$. The map $f$ has
  the homotopy LLP with respect to
  $$
  g^*:\map_*(D,Z)\to \map_*(C,Z)
  $$
  if and only if $g$ has the homotopy LLP with respect to
  $$
  f^*:\map_*(B,Z)\to \map_*(A,Z)
  $$
\end{lem}
\begin{pf}
  We use adjointness to note that
  the existence of a dashed lift in the diagram
  $$
  \xymatrix{
    A \ar[d]^f \ar[r] & {\map_*(D, Z)} \ar[d]^{g^*} \\
    B \ar[r] \ar@{-->}[ur] & {\map_*(C, Z)}
    }
  $$
  is equivalent to the existence of the dashed map in the following
  diagram.
  $$
  \xymatrix{
    {A\wedge D} \ar[dd]_{f\wedge id} \ar[dr] &&
    {A\wedge C} \ar[ll]_{id\wedge g}
    \ar[dd]^{f\wedge id} \\
    & Z \\
    {B\wedge D} \ar@{-->}[ur] && {B\wedge C} \ar[ul]
    \ar[ll]_{id\wedge g}
    }
  $$
  This in turn is equivalent to the lifting property as indicated
  on the next diagram.
  $$
  \xymatrix{
    C \ar[d]^g \ar[r] & {\map_*(B, Z)} \ar[d]^{f^*} \\
    D \ar[r] \ar@{-->}[ur] & {\map_*(A, Z)}
    }
  $$
\qed \end{pf}

\begin{lem}  \label{lemma_llp}
  Let $g:\bigvee_{n\geq 0} S^n \to \bigvee_{n\geq 0} S^n $
  be the trivial map. A map $f:A \to B$ is a $Z$-equivalence
  if and only if it has the homotopy LLP with respect to
  $$
  g_+^*:\map_*((\bigvee_{n\geq 0} S^n)_+,Z) \to \map_*((\bigvee_{n\geq 0} S^n)_+,Z)
  $$
\end{lem}
\begin{pf}
  By Lemma \ref{lemma_symmetry} $f$ has the homotopy LLP with respect
  to $g_+^*$ if and only if $g_+$ has the homotopy LLP with respect to
  $f^*:\map_*(B,Z)\to \map_*(A,Z)$.
  Obviously if $f^*$ is a weak equivalence then $g_+$ has the
  homotopy LLP hence the proof will be complete once we show that
  the homotopy LLP for $g_+$ implies that $f^*$ is a weak equivalence.
  We see that if $g_+$ has the homotopy LLP with respect to $f^*$
  then all the maps $g^n_+:S^n_+\to\{*\}_+\to S^n_+$ for $n\geq 0$
  have the homotopy LLP.
  The case $n=0$ implies that $f^*$ induces a bijection on the components.

  We are proving that $f^*$ induces isomorphisms of homotopy groups
  of the corresponding components.
  Assume that $f$ is an inclusion $A\hookrightarrow B$ of simplicial
  sets and $Z$ is a fibrant simplicial set. We fix any map
  $b_0:B\to Z$ as a basepoint of $\map_*(B,Z)$ and $a_0=f^*(b_0)$
  as a basepoint of $\map_*(A,Z)$.
  The homotopy LLP for $g^n_+$ for $n>0$ implies that
  $f^*$ induces bijections of the homotopy groups modulo the action
  of the fundamental group:
  $$
  \pi_n(\map_*(B,Z),b_0)/\!\!\sim
    \;\to \pi_n(\map_*(A,Z),a_0)/\!\!\sim
  $$
  Since $0$ is fixed by the action of the fundamental group
  we see that
  $$
  f^*_n:\pi_n(\map_*(B,Z),b_0) \hookrightarrow \pi_n(\map_*(A,Z),a_0)
  $$
  is a monomorphism for $n>0$.
  Choose an element $\tilde{\alpha}\in \pi_n(\map_*(A,Z),a_0)$.
  It is represented by some $\alpha:A\wedge S^n_+\to Z$
  such that $\alpha|{\empty\atop {A = A\wedge\{*\}_+}} = a_0$.
  We construct the following diagram.
  $$
  \xymatrix{
    A\wedge S^n_+ \ar[dr]^\alpha \ar@{_{(}->}[dd]_{f\wedge id} &&
      {A\wedge S_+^n} \ar[ll]_{id\wedge g^n_+}
                      \ar@{_{(}->}[dd]^{f\wedge id} \\
    & Z \\
    B\wedge S^n_+ \ar@{-->}[ur]^\beta && B\wedge S^n_+ \ar[ul]_b \ar[ll] \\
  }
  $$
  The map $b$ is the composition
  $B\wedge S^n_+ \to B\wedge \{*\}_+ = B \stackrel{b_0}{\to} Z$.
  The diagram commutes by the definition of $a_0$ as $b_0f$.
  By the proof of Lemma \ref{lemma_symmetry} the assumption that
  $g^n_+$ has the homotopy LLP with respect to $f^*$ implies
  the existence of the dashed map $\beta$ which closes
  this diagram up to homotopy.
  Since $f^*$ is a bijection on components we see that
  $\beta|{\empty\atop B\wedge\{*\}_+}:B\to Z$ must be homotopic to $b_0$.
  Since $A\wedge\{*\}_+\hookrightarrow B\wedge S^n_+$ is a cofibration
  we can find $\beta_1$, homotopic to $\beta$, such that
  $\beta_1|{\empty\atop B\wedge\{*\}_+} = b_0$.
  We see that $\beta_1$ induces an element $\tilde{\beta}$ in
  $\pi_n(\map_*(B,Z),b_0)$ such that $f^*(\tilde{\beta})=\tilde{\alpha}$
  hence $f^*$ is a weak equivalence.
\qed \end{pf}

\section{Categories of pairs and topologized objects}
\label{sec_categories}

In this section we collect some categorical definitions and facts
which will be used in section \ref{section_at_a_space}.
Some statements refer to a general category $\mathcal{C}$,
however for us the interesting
cases are when $\mathcal{C}=\mathcal{S}_*$ (pointed simplicial sets)
or $\mathcal{C}=Ho_*$.

\begin{defn}  \label{definition_C2}
  Given a category $\mathcal{C}$ we will denote by
  $\mathcal{C}^2$ the usual {\em category of pairs} whose objects
  are the maps in $\mathcal{C}$ and whose maps are commutative
  squares in $\mathcal{C}$ as below.
  $$
  \xymatrix{
    A \ar[d]_f \ar[r]^{h_A} & S \ar[d]^g \\
    B \ar[r]^{h_B} & T
  }
  $$
\end{defn}

Following Bousfield and Friedlander
(see \cite{bousfield_friedlander} A3) we introduce
a model category structure on $\mathcal{C}^2$.

\begin{defn}
  Let $\mathcal{C}$ be a model category. A map $h:f\to g$ as in
  Definition \ref{definition_C2} is called a {\em weak equivalence}
  (respectively {\em fibration}) if both $h_A$ and $h_B$ are
  weak equivalences (respectively fibrations). It is a
  {\em cofibration} if $h_A:A\to S$ and
  $(h_B,g):B\amalg_A S \to T$ are cofibrations.
  This implies that $h_B:B \to T$ is also a cofibration.
\end{defn}
Note that an object $f:A\to B$ is cofibrant in $\mathcal{C}^2$
if $A$ is cofibrant in $\mathcal{C}$ and the map $f$ is
a cofibration in $\mathcal{C}$. It is fibrant if both
$S$ and $T$ are fibrant in $\mathcal{C}$.

We will be interested in $ho\mathcal{S}^2_*$ the homotopy
category of pairs when $\mathcal{C}=\mathcal{S}_*$ the category of
pointed simplicial sets.
The obvious functor $F:ho\mathcal{S}^2_*\to Ho^2_*$
induces equivalence of categories.

\comment{
\begin{lem}
  The obvious functor $F:ho\mathcal{S}^2_*\to Ho^2_*$
  induces equivalence of categories.
\end{lem}
\begin{pf}
  We define the ``inverse'' functor
  $D:Ho^2_* \to ho\mathcal{S}^2_*$ as follows.
  Given an object $f:[A]\to [B]$ in $Ho^2_*$ we represent it
  by a map $\overline{f}:\overline{A}\to\overline{B}$ in $\mathcal{S}_*$
  with $\overline{A}$ and $\overline{B}$ fibrant.
  This way $\overline{f}$ is a fibrant object in $\mathcal{S}^2_*$.
  Let $D(f)$ be the class represented by $\overline{f}$.
  Obviously the composition $FD$ is an identity on $Ho^2_*$.

  We need a natural equivalence between $Id_{ho\mathcal{S}^2_*}$
  and $DF$. Consider a cofibrant in $\mathcal{S}^2_*$ object
  $f:A\hookrightarrow B$. We have $F(f):[A]\to[B]$. Let $DF(f)$
  be represented by $\overline{f}:\overline{A}\to\overline{B}$.
  Since $F(f)=FDF(f)$ and $A$ and $B$ are cofibrant and
  $\overline{A}$ and $\overline{B}$ are fibrant we obtain
  a homotopy commutative diagram
  $$
  \xymatrix{
    A \ar@{_{(}->}[d]_f \ar[r]^{h_A} & \overline{A} \ar[d]^{\overline{f}} \\
    B \ar[r]^{h_B} & \overline{B}
  }
  $$
  where $h_A$ and $h_B$ induce identity maps in $Ho_*$. Since $f$
  is a cofibration and $\overline{B}$ is fibrant we can modify
  $h_B$ within its homotopy class so that the diagram commutes in
  $\mathcal{S}^2_*$. This way we obtain a map $h:f\to\overline{f}$
  in $\mathcal{S}^2_*$ which is an equivalence in $\mathcal{S}^2_*$.
\qed \end{pf}
}

Some of the definitions below are chosen after \cite{deleanu}.
For any category $\mathcal{C}$ and an object $X$ of $\mathcal{C}$
a {\em topologized object over} $X$ is a factorization
$$
\xymatrix{
  & {\mbox{Top}} \ar[dr]^G \\
  {\mathcal{C}^{op}} \ar[ur]^{X^{\#}} \ar[rr]_{\mathcal{C}(-,X)}
    && {\mbox{Sets}} \\
}
$$
where $G$ is the forgetful functor. We say that a morphism
$f:X\to Y$ is {\em continuous} if it induces a natural transformation
$f^{\#}:X^{\#}\to Y^{\#}$, that is to say, the map
$\hom_{\mathcal{C}}(Z,f)$ is continuous with respect to the
topologies of $X^{\#}Z$ and $Y^{\#}Z$ for all $Z$ in $\mathcal{C}$.
\begin{defn}\label{definition_compact}
  We say that a topologized object $X$ is {\em compact} if
  the corresponding functor $X^{\#}$ takes values in
  compact Hausdorff spaces. A category of compact objects
  and continuous morphisms in $\mathcal{C}$ will be denoted
  by $C\mathcal{C}$
\end{defn}

\begin{lem}  \label{lemma_compact_map}
  If $g:S\to T$ is a map in $CHo_*$ then it is naturally a compact
  object in $Ho^2_*$. In other words the categories
  $(CHo_*)^2$ and $CHo^2_*$ have the same objects.
\end{lem}
\begin{pf}
  We need to show that for any $f:A\to B$ in $Ho_*$ the set
  $\hom_{Ho^2_*}(f,g)$ has a natural compact topology.
  This is obvious since this set is the limit of the
  following diagram
  $$
  \xymatrix{
    {[A,S]\times [B,T]}
    \ar@<0.5ex>[r]^(0.65)\varphi \ar@<-0.5ex>[r]_(0.65)\psi &
    {[A,T]}
  }
  $$
  where the entries are compact since $S$ and $T$ are in $CHo_*$.
  The maps $\varphi(\alpha,\beta)=g\alpha$ and
  $\psi(\alpha,\beta)=\beta f$ are continuous.
\qed \end{pf}

By adjointness argument we immediately obtain the following.

\begin{lem}  \label{lemma_compact_mapping_complex}
  If $T$ is in $CHo_*$ then for any $X$ the space $\map_*(X,T)$
  is in $CHo_*$ and for any map $f:X\to Y$ the induced map\
  $\map_*(Y,T)\to \map_*(X,T)$ is continuous.
\end{lem}

\section{Localizations at a space}  \label{section_at_a_space}

In this section we will prove
(Theorem \ref{theorem_localization_at_Z})
that localization at a space $Z$ exists whenever $Z$
is a homotopy retract of
a compact object in the sense of Definition \ref{definition_compact}.
We attain this by showing that
for such spaces $Z$ any $Z$-equivalence can be presented as a filtered
colimit of $Z$-equivalences of bounded cardinalities so that we can
use Lemma \ref{lemma_set_generated_equivalences}.

Let $\mathcal{S}^2_*$
be the usual category of maps in $\mathcal{S}_*$.
We will say that $f_0$ is a subobject of $f$ if there is a
cofibration $f_0\hookrightarrow f$ and will denote this fact
by $f_0\subseteq f$.
Given $f:A\to B$ we will write $|f|$ for the number
of nondegenerate simplexes of $A\vee B$ and will say that $f$ is
finite if $|f|$ is.

\begin{lem}  \label{lemma_Z_extensions}
  Let $f\subseteq h$ be cofibrant objects in $\mathcal{S}^2_*$.
  Let $g$, fibrant in $\mathcal{S}^2_*$, represent an object in
  $Cho\mathcal{S}^2_*$.
  Let $\alpha \in \hom_{\mathcal{S}^2_*}(f,g)$. If for every
  finite subobject $k\subseteq h$ the map $\alpha$
  extends to $f\cup k$
  then $\alpha$ extends to $h$.
\end{lem}
\begin{pf}
  Let $t$ be in $\mathcal{S}^2_*$
  such that $f\subseteq t\subseteq h$.
  Let $\nolinebreak{r:\hom_{ho\mathcal{S}^2_*}(t,g)\to
  \hom_{ho\mathcal{S}^2_*}(f,g)}$
  be the restriction map. Define $E(t)$
  as $r^{-1}([\alpha])$ that is the set of all extensions,
  in $ho\mathcal{S}^2_*$, of $\alpha$ to
  $t$. Since $r$ is a continuous map between compact spaces we see
  that $E(t)$ is empty or compact.
  The limit $\lim E(f\cup k)$ taken over all finite subobjects of
  $h$ is nonempty since it is directed and the sets $E(f\cup k)$
  are compact (nonempty by assumption).
  The proof will be complete once we show that $E(h)$ is nonempty.
  We will show that $E(h)=\lim E(f\cup k)$.
  Let $\map_*(t,g)$ be a simplicial set whose $n$-simplexes
  form a set $\hom_{\mathcal{S}^2_*}(t\wedge(\Delta^n_+),g)$
  and whose faces and degeneracies are induced by the cosimplicial
  structure on $\Delta^{\stackrel{\bullet}{\empty}}$.
  Obviously $\pi_0(\map_*(t,g))=E(t)$.
  Since $g$ represents an object in $Cho\mathcal{S}^2_*$ we see that
  $\pi_q(\map_*(t,g))=\hom_{ho\mathcal{S}^2_*}
  (t\wedge(\Delta^q/\partial\Delta^1),g)$
  is compact for $q\geq 0$ which gives us the last equation in
  the following sequence.
  $$
  \pi_0(\map_*(h,g))
  =\pi_0(\map_*(\colim f\cup k,g))
  =\pi_0(\map_*(\hocolim f\cup k,g))=
  $$
  $$
  =\pi_0(\holim\;\map_*(f\cup k,g))
  =\lim \pi_0(\map_*(f\cup k,g))
  $$
  This means that
  $$
  E(h)=\lim E(f\cup k)
  $$
\qed \end{pf}

Directly from Lemma \ref{lemma_Z_extensions} we obtain the following
statement.

\begin{lem}\label{lemma_Z_obstruction}
  Given cofibrant $f$ and fibrant $g$ in $\mathcal{S}^2_*$ with $g$
  representing an object in $Cho\mathcal{S}^2_*$
  there is a cardinal number $\tau=\tau(f,g)$ such that for
  any $h$ in $\mathcal{S}^2_*$ with $f\subseteq h$
  there is $k$ in $\mathcal{S}^2_*$ such that
  $f\subseteq k\subseteq h$ and $|k|\leq\tau$ and
  if $\alpha:f\to g$ extends to $\alpha_k:k\to g$ then
  it extends to $\alpha_h:h\to g$.
\end{lem}

\begin{pf}
  For each $\alpha:f\to g$ which does not factor as
  $f\hookrightarrow h \to g$ Lemma \ref{lemma_Z_extensions} gives us a
  finite object $k_\alpha$ in $\mathcal{S}^2_*$
  such that $\alpha$ does not factor as
  $f\hookrightarrow f\cup k_{\alpha} \to g$. We can take
  $k=f\cup\bigcup_{\alpha} k_{\alpha}$.
  Since each $k_\alpha$ is finite and the number of possible
  maps $\alpha$ depends only on $f$ and $g$ we see that there is an
  upper bound for the cardinality of $k$ which depends only on $f$ and
  $g$.
\qed \end{pf}

The role of this Lemma is following. We think of $f$ and $g$ as
fixed and of $h$ as uncontrollably big. We want the obstruction to
extending a map from $f$ to $h$ to be detected on some $k$ whose
cardinality we can control.

\begin{lem}\label{lemma_Z_onto}
  Given cofibrant $f$ and fibrant $g$ in $\mathcal{S}^2_*$ with $g$
  representing an object in $Cho\mathcal{S}^2_*$
  there is a cardinal number $\delta=\delta(f,g)$ such that for
  any $h$ in $\mathcal{S}^2_*$ with $f\subseteq h$
  there is $k$ in $\mathcal{S}^2_*$ such that
  $f\subseteq k\subseteq h$ and $|k|\leq\delta$ and
  the restriction map
  $\hom_{ho\mathcal{S}^2_*}(h,g)\twoheadrightarrow
  \hom_{ho\mathcal{S}^2_*}(k,g)$
  is an epimorphism.
\end{lem}

\begin{pf}
  The object $k$ is constructed as a union of an ascending
  chain $f=k_0\subseteq k_1\subseteq ...
  \subseteq k_n \subseteq ... $. This chain is built by induction on
  $n$. Given $k_n$ we use Lemma \ref{lemma_Z_obstruction} to choose
  $k_{n+1}$ so that $k_n \subseteq k_{n+1} \subseteq h$
  and if a map $k_n\to g$ extends to $k_{n+1}$ then it extends to
  $h$.

  Given $\alpha:k\to g$ we need to show that we can extend $\alpha$ to
  $\tilde{\alpha}:h \to g$. By the construction of $k$ there are maps
  $\alpha_n:h\to g$ such that
  $\alpha_n|{\empty\atop k_n} \simeq \alpha|{\empty\atop k_n}$.
  Since by assumption
  $\hom_{ho\mathcal{S}^2_*}(h,g)$ is compact
  we can take $\tilde{\alpha}$ to be an accumulation point of the set
  $\{\alpha_n\}$.

  We have
  $\tilde{\alpha}|{\empty \atop k_n} \simeq \alpha|{\empty \atop k_n}$
  for all $n$ since the sequence
  $\alpha_i|{\empty \atop k_n}\in\hom_{ho\mathcal{S}^2_*}(k_n,g)$
  converges to $\alpha|{\empty \atop k_n}$, it is actually constant
  for $i\geq n$, and the restriction map
  $\hom_{ho\mathcal{S}^2_*}(h,g)\to\hom_{ho\mathcal{S}^2_*}(k_n,g)$
  is continuous.

  A similar argument as in the last paragraph of the proof
  of Lemma \ref{lemma_Z_extensions} tells us that
  $$
  \alpha \in \hom_{ho\mathcal{S}^2_*}(k,g) =
  \lim \hom_{ho\mathcal{S}^2_*}(k_n,g)
  $$
  hence $\tilde{\alpha}|{\empty \atop k_n} \simeq \alpha|{\empty \atop k_n}$
  for all $n$ implies
  $\tilde{\alpha}|{\empty \atop k} \simeq \alpha$.
\qed \end{pf}

\begin{lem}\label{lemma_colimit_of_equivalences}
  Let $g$ in $\mathcal{S}^2_*$ represent an object in
  $Cho\mathcal{S}^2_*$. Let cofibrant $h$ and fibrant $p$ be in
  $\mathcal{S}^2_*$. Let $p$ be a retract in $\mathcal{S}^2_*$
  of $g$ and $h$ have the homotopy LLP with respect to $p$.
  There is a cardinal $\gamma=\gamma(g)$ such that $h$ is a colimit of
  subobjects $h_\alpha$ such that each $h_\alpha$ has the homotopy LLP
  with respect to $p$ and $|h_\alpha|\leq \gamma$.
\end{lem}
\begin{pf}
  We can write $h$ as $h=\colim h_\alpha$ where each $h_\alpha$ is
  finite. Inductively we replace $h_\alpha$ with objects $h_{*\alpha}$
  that have the left lifting property with respect to $p$. We start with
  the trivial object in $\mathcal{S}^2_*$, a map between
  spaces consisting of a basepoint only, which need not be
  replaced. Suppose that for some $\alpha_0$ all subobjects of
  $h_{\alpha_0}$ have been replaced. Let
  $h'=h_{\alpha_0}\cup\bigcup_{\alpha<{\alpha_0}}h_{*\alpha}$.
  Lemma \ref{lemma_Z_onto} gives us a factorization
  $$
  h' \hookrightarrow h_{*\alpha_0} \hookrightarrow h
  $$
  such that the restriction map

  \begin{equation}\label{equation_hom_onto}
  \hom_{ho\mathcal{S}^2_*}(h,g)\twoheadrightarrow
  \hom_{ho\mathcal{S}^2_*}(h_{*\alpha_0},g)
  \end{equation}

  is an epimorphism.
  We want to show that $h_{*\alpha_0}$ has the homotopy LLP
  with respect to $p$.
  For any map $\varphi:h_{*\alpha_0}\to p$ consider a diagram
  $$
  \xymatrix{
    {h_{*\alpha_0}} \ar@{^{(}->}[d] \ar[r]^\varphi & p \ar[d] \\
    h \ar[r]^\psi & g \ar@/_1pc/[u]
    }
  $$
  where the map $\psi$ exists by (\ref{equation_hom_onto}). Since by
  assumption $h$ has the left lifting property with respect to $p$
  and any map from $h_{*\alpha_0}$ to $p$ factors through $h$ we
  obtain the homotopy LLP for $h_{*\alpha_0}$
  with respect to $p$.
  We see that $|h_{*\alpha_0}|$ depends only
  on $g$, on $h_{*\alpha}$ for $\alpha<\alpha_0$ and on the bounds
  $\delta(h_{*\alpha},g)$ from Lemma \ref{lemma_Z_onto}.
\qed \end{pf}

We are ready to prove the main theorem of this paper.
In the following we prefer to work in the $Ho^2_*$
rather than in the equivalent category $ho\mathcal{S}^2_*$.
\begin{thm}\label{theorem_localization_at_Z}
  Let $\overline{Z}$ in $Ho_*$ represent an object in $CHo_*$.
  For any $Z$ in $Ho_*$, a homotopy retract of $\overline{Z}$,
  there exists a
  map $f$ such that $L_f$ is a localization at $Z$.
\end{thm}
\begin{pf}
  To use Lemma \ref{lemma_llp} we consider maps
  $$
  p:\map_*((\bigvee_{n\geq 0} S^n)_+,Z) \to \map_*((\bigvee_{n\geq 0} S^n)_+,Z)
  $$
  and
  $$
  g:\map_*((\bigvee_{n\geq 0} S^n)_+,\overline{Z}) \to
  \map_*((\bigvee_{n\geq 0} S^n)_+,\overline{Z}).
  $$
  We observe that $p$ is a homotopy retract of $g$ and by
  Lemma \ref{lemma_compact_mapping_complex} $g$ represents an object
  in $CHo^2_*$.
  By Lemma
  \ref{lemma_llp} a map $h$ is a $Z$-equivalence
  if and only if it has the homotopy LLP
  with respect to $p$. By Lemma
  \ref{lemma_colimit_of_equivalences}
  there is a cardinal $\gamma=\gamma(g)$ such that
  any $Z$-equivalence $h$ is a colimit of
  $Z$-equivalences whose
  cardinalities do not exceed $\gamma$. Since this is a directed
  colimit of cofibrations it is equivalent to a homotopy colimit. By
  Lemma \ref{lemma_set_generated_equivalences} we can take $f$ to be a
  wedge of all $Z$-equivalences whose cardinality does not exceed
  $\gamma$.
\qed \end{pf}

Since one would like to remove the compactness assumption in
Theorem \ref{theorem_localization_at_Z} we briefly review the points
where we used it in the proof.
The key property we used in Lemmas \ref{lemma_Z_extensions} and
\ref{lemma_Z_onto} is that for a compactly topologized $C$ and a
directed diagram $X_i$ in $Ho_*$ there is a bijection
$$
[\holim X_i, C] \stackrel{\simeq}{\to} \lim [X_i,C]
$$
Other properties are much simpler, in Lemma \ref{lemma_Z_onto}
we needed to know that an infinite subset of a compact topological space
has an accumulation point and in Lemma \ref{lemma_compact_map}
that a closed subspace of a product of compact spaces is compact.

We end this section with Example \ref{example_noncompact_retract}
which shows that the "retract" condition in Theorem
\ref{theorem_localization_at_Z} is relevant. More precisely
there are spaces which represent objects
in $CHo_*$ but whose retracts are not in $CHo_*$.

We will need the following two lemmas.
By a {\em simplicial compact space} we understand a simplicial object
in the category of compact (Hausdorff) topological spaces.

\begin{lem}  \label{lemma_first_topology}
  Let $X$ be a simplicial set and $Z$ a simplicial compact space.
  The set $\hom_{\mathcal{S}_*}(X,Z)$ has a natural compact topology.
\end{lem}
\begin{pf}
  To see this observe that $\hom_{\mathcal{S}_*}(X,Z)$
  is a subset of
  $$
  \displaystyle\prod_n\mbox{Sets}(X_n,Z_n)\cong
  \prod_n\prod_{X_n}Z_n
  $$
  which has a compact product topology. The subset
  $\hom_{\mathcal{S}_*}(X,Z)$
  is determined by a number of equations (see May
  \cite{may} 1.2)
  between continuous maps so it forms a closed hence compact subspace
  of the product.
\qed \end{pf}

\begin{lem}  \label{lemma_compact_mapping_space}
  Let $T$ be a simplicial compact space which is fibrant
  as a simplicial set. Then $T$ naturally represents an object
  in $CHo_*$.
\end{lem}

\begin{pf}
  We need to show that for any simplicial set $X$ the set
  $[X,T]$ is naturally compact. We have
  $\map_*(X,T)_k = \hom_{\mathcal{S}_*}(X\wedge(\Delta^k_+),T)$
  hence by Lemma \ref{lemma_first_topology} the mapping space
  $\map_*(X,T)$ is a simplicial compact space.
  Since $[X,T]=\pi_0\map_*(X,T)$ hence by
  Proposition $4.7$ in \cite{bousfield_israel}
  it is naturally compact.
\qed \end{pf}

\begin{example}  \label{example_noncompact_retract}
  Let $n>0$, $Z=K(\mathbb{Q},n)$ and $\overline{Z}=K(S^1,n)$.
  As a model of $K(S^1,n)$ we use the one described in $1.2$ of
  \cite{bousfield_unpublished}; $K(S^1,n)_t$ is a product of
  $({t\atop n})$ copies of $S^1$, hence it is a compact
  topological space, faces and degeneracies are given
  by projections and group operations hence they are continuous.
  This model of $K(S^1,n)$ is a simplicial compact space
  which is fibrant as a simplicial set.
  It has a homotopy type of an Eilenberg-Mac Lane space for
  $S^1$ viewed as a discrete group.
  The group $S^1$ is a direct sum of
  $\mathbb{Q}/\mathbb{Z}$ and a rational vector space hence
  $\mathbb{Q}$ is a retract of $S^1$ and so $Z$ is a retract of
  $\overline{Z}$.
  We have $\overline{Z}$ which represents an object
  in $CHo_*$ and its retract $Z$ which does not represent any
  objects in $CHo_*$
  since $\pi_nZ=\mathbb{Q}$ is an infinite countable group
  hence admits no compact structure.
\end{example}

\section{Applications and Examples}\label{section_applications_and_examples}

We note that Theorem \ref{theorem_localization_at_Z} implies
the existence of localizations at spaces which belong
to the following classes:
\begin{itemize}
  \item[a)]
    Profinite completions of other spaces.
  \item[b)]
    Simplicial compact spaces which are fibrant as simplicial sets
    (Lemma \ref{lemma_compact_mapping_space}).
  \item[c)]
    Mapping spaces with targets in a) or b)
    (Lemma \ref{lemma_compact_mapping_complex}).
\end{itemize}

Our first example of a localization at a space
is an idempotent approximation to
the profinite completion. The work of Rao \cite{rao}
implies the existence of such an approximation defined on the
nilpotent spaces. Here we don't require such assumptions.

The profinite completion was introduced by Sullivan in section $3$
of \cite{sullivan} via the Brown representability theorem.
To a given space $X$ he assigns another space $\hat{X}$
which represents the functor
$\hat{X}(Y)=\lim_{(X\downarrow\mathcal{F})}[Y,F]$.
The limit is taken over the category
$(X\downarrow \mathcal{F})$ whose objects are maps $X\to F$
in $Ho_*$ with $F$ connected
and $\pi_qF$ finite for all $q>0$.
The morphisms are
commutative diagrams in $Ho_*$ as below.
$$
\xymatrix{
  & X \ar[dl] \ar[dr] \\
  F_1 \ar[rr] && F_2 \\
}
$$
The functor $F:(X\downarrow \mathcal{F})\to \mathcal{S}_*$
takes an object $X\to F_0$ to the space $F_0$.
This limit is
well defined since the category
$(X\downarrow \mathcal{F})$ is equivalent to a small category.

\begin{thm}\label{approximation_to_completion}
  There exists an idempotent approximation to the profinite
  completion. More precisely, there is the terminal localization among
  localizations $L$ which admit the following factorization.
  $$
  X\to LX\to \hat{X}
  $$
\end{thm}
\begin{pf}
  For each homotopy class of connected spaces with
  $\pi_qF$ finite for all $q>0$ choose a representative $F$.
  Let $Z=\prod F$ be the product of those representatives.
  Since each $F$ is naturally compact (in the sense of
  Definition \ref{definition_compact}) and
  $[Y,Z]=\prod[Y,F]$ for all $Y$ we see that $Z$ is compact.
  The localization $L_Z$ exists by Theorem
  \ref{theorem_localization_at_Z}. We observe that if $F$ is
  connected with $\pi_qF$ finite for $q>0$ then $F$ is
  $Z$-local. Let $r:Z\to F \hookrightarrow Z$ be the retraction onto
  the axis that corresponds to $F$. We see that
  $F\simeq \holim(...\stackrel{r}{\to} Z\stackrel{r}{\to} Z)$ hence
  it is $Z$-local.
  This implies that
  $[L_ZX,F]\to[X,F]$ is a bijection and consequently
  that the categories $(X\downarrow \mathcal{F})$ and
  $(L_ZX\downarrow \mathcal{F})$ are equivalent hence
  $\hat{X}\simeq (L_ZX)\hat{}$ which leads us to the
  factorization we were looking for:
  \begin{equation}\label{equation_profinite_factorization}
  X\to L_ZX\to (L_ZX)\hat{}\simeq\hat{X}
  \end{equation}

  It remains to show that $L_Z$ is the terminal localization which
  admits factorization (\ref{equation_profinite_factorization}).
  Suppose that a localization $T$ also admits
  (\ref{equation_profinite_factorization}). Since profinite completion
  is idempotent on finite spaces $F$ as above we have
  $$
  F \to TF \to \hat{F} \simeq F
  $$
  so $F$ is a homotopy retract of $TF$ hence $T$-local. This means
  that the space $Z$ is $T$-local hence by the definition of $L_Z$ we
  have $T\leq L_Z$.
\qed \end{pf}

\begin{thm}\label{unstable_cohomological}
  Let $h^*$ be a cohomology theory represented by an $\Omega$-spectrum
  $\{\underline{h}_n\}$.
  If each $\underline{h}_n$ is a homotopy
  retract of a compact, in the sense of Definition
  \ref{definition_compact}, space then
  there exists a map $f$ such that $L_f$-equivalences and
  $h^*$-equivalences coincide. In particular the corresponding
  cohomological localization exists.
\end{thm}
\begin{pf}
  Let $Z=\prod \underline{h}_n$ and use
  Theorem \ref{theorem_localization_at_Z}.
\qed \end{pf}

\end{document}